\documentclass[11pt]{article}

\usepackage{amsmath}
\usepackage{amssymb}
\usepackage{amsbsy}
\usepackage{amsfonts}
\usepackage{mathtools}
\usepackage{bbm}
\usepackage{comment}
\usepackage{color}
\newtheorem{theorem}{Theorem}[section]%
\newtheorem{lemma}[theorem]{Lemma}%
\newtheorem{prop}[theorem]{Proposition}%

\newenvironment{pf}{\medskip\noindent{Proof:}
  \hspace{-.5cm}      \enspace}{\hfill \qed \newline \smallskip}

   \def\z{\zeta}

\def\f{\noindent}

\def\mod{\hbox{\rm mod}\;}

\newcommand{\qed}{\mbox{\raisebox{0.7ex}{\fbox{}}} \vspace{4truemm}}

\begin{document}

\begin{center}
\Large{\textbf{On the complexity of Cayley graphs on a dihedral group}}
\end{center}

\vskip 5mm {\small
\begin{center}

{\textbf{Bobo Hua,}}\footnote{Fudan University, bobohua@fudan.edu.cn}
{\textbf{A.~D.~Mednykh,}}\footnote{{\small\em Sobolev Institute of Mathematics, Novosibirsk State University, smedn@mail.ru}} 
{\textbf{I.~A.~Mednykh,}}\footnote{{\small\em Sobolev Institute of Mathematics, Novosibirsk State University, ilyamednykh@mail.ru}}
{\textbf{Lili Wang,}}\footnote{School of Mathematics and Statistics \& Key Laboratory of Analytical Mathematics and Applications (Ministry of Education) \& Fujian Provincial Key Laboratory of Statistics and Artificial Intelligence \& Fujian Key Laboratory of Analytical Mathematics and Applications (FJKLAMA) \& Center for Applied Mathematics of Fujian Province (FJNU), Fujian Normal University, liliwang@fjnu.edu.cn }
\end{center}}

\title{ \vspace{-1.2cm}
$On the complexity of Cayley  graph on a dihedral group.$
\thanks{Supported by }}

\begin{abstract}
In this paper, we investigate the complexity of an infinite family of Cayley graphs $\mathcal{D}_{n}=Cay(\mathbb{D}_{n}, b^{\pm\beta_1},b^{\pm\beta_2},\ldots,b^{\pm\beta_s}, a b^{\gamma_1}, a b^{\gamma_2},\ldots, a b^{\gamma_t} )$ on the dihedral group $\mathbb{D}_{n}=\langle a,b| a^2=1, b^n=1,(a\,b)^2=1\rangle$ of order $2n.$ 

We obtain a closed formula for the number $\tau(n)$ of spanning trees in $\mathcal{D}_{n}$ in terms of Chebyshev polynomials, investigate  some arithmetical  properties of this function,  and find its asymptotics as $n\to\infty.$  Moreover, we show that the generating function $F(x)=\sum\limits_{n=1}^\infty\tau(n)x^n$  is a rational function with integer coefficients.
\bigskip

\f\textbf{AMS classification:} 05C30, 39A12
\end{abstract}

\section*{Introduction}

Let $G$ be a finite connected graph. The notion of the complexity of a graph can be defined in several different ways. One can consider the number of edges or vertices, the number of spanning trees or rooted spanning forests. The important characteristic of a graph is the Kirchhoff index defined as the sum of resistance distances between vertices. All the above-mentioned values can be expressed in terms of the Laplacian spectrum of a graph. In particular, by the famous Kirchhoff Matrix-Tree Theorem the number of spanning trees in a connected graph is equal to the product of all non-zero eigenvalues of its Laplacian matrix divided by the number of vertices.         

The study of such invariants usually leads to the following question: how to find the product of eigenvalues of the Laplacian matrix? If the size (number of vertices) of a graph is small, it is an easy task. However, the most interesting cases involve the family of graphs with increasing number of vertices. The direct calculation of this product becomes tedious and unmanageable when the number of vertices $n$ of the graph tends to infinity. To solve this problem, we use the techniques developed in previous papers by the authors \cite{GrunKwonMed}, \cite{KwonMedMed} and \cite{Med1}. 

As a result, one can find a closed formula which is the product of a bounded number of factors, each given by the $n$-th Chebyshev polynomial of the first kind evaluated at the roots of some polynomial of prescribed degree. This paves the way to investigate arithmetical properties and asymptotics.   

The complexity of a graph plays an important role in statistic physics, where the graphs with arbitrarily large number of vertices are considered (\cite{JacSalSok}, \cite{SW00}, \cite{Wu77}). With increasing number of vertices the structure of the Laplacian characteristic polynomial becomes quite complicated. In this case, the most interesting invariants are given by their asymptotics. See surveys for counting spanning trees \cite{BoePro}, \cite{MedMedSurvey} and the references therein. 

The aim of the present paper is to produce explicit analytic formulas for the number of spanning trees in a Cayley graph $\mathcal{D}_n,$ see \eqref{eq:d1} for the definition, on a dihedral group. Such formulas help in the investigation of number--theoretical properties and asymptotics of several spectral invariants of the graph. This research expands series of publications by various authors on the complexity of circulant graphs (\cite{AbrBaiMed}, \cite{GrunMedArs22}, \cite{GrunKwonMed}, \cite{JustLou15}, \cite{MedMed2018}, \cite{MedMedK}, \cite{XiebinLinZhang}, \cite{ZhangYongGol}, \cite{ZhangYongGolin}). Note that the circulant graph is a Cayley graph on a cyclic group.  

The methods laid out in the present paper can be equally used to find explicit formulas for the number of rooted spanning forests and Kirchhoff index.  

The paper is organized as follows. Some preliminary results and basic definitions are given in Section~\ref{basic}. Section~\ref{spectra} is devoted to the investigation of spectral properties of the Laplacian matrix of the Cayley graph $\mathcal{D}_n$ on the dihedral group $\mathbb{D}_n.$ In Section~\ref{counting}, we present explicit formulas for the number $\tau(n)$ of spanning trees in graph $\mathcal{D}_n.$ The formulas will be given in terms of Chebyshev polynomials evaluated in the roots of Laurent polynomial associated with $\mathcal{D}_n.$ In Section~\ref{arithmetic}, we provide some arithmetical properties of function $\tau(n)$ for the family $\mathcal{D}_n.$ More precisely, we show that the number of spanning trees  can be represented in the form $\tau(n)=p\,n\,a(n)^2,$ where $a(n)$ is an integer sequence and $p$ is a prescribed number depending only on parameters of $\mathcal{D}_n$ and the parity of $n.$ In Section~\ref{asymptotic}, we use explicit formulas for the complexity to produce its asymptotic. Section~\ref{genfuction} describes the structure of the generation function $F(x)=\sum_{n=1}^{\infty}\tau(n)x^{n}.$ We show that $F(x)$ is a rational function with integer coefficients satisfying some symmetry property. In the last section, we illustrate the obtained results by a series of examples.

\section{Basic definitions and preliminary facts}\label{basic}

Consider a finite connected graph $G$ with possibly multiple edges, but without loops. We denote by $V(G)$ and $E(G)$ the set of vertices and the set of edges of $G$ respectively. 

A \textit{tree} is a connected undirected graph without cycles. A \textit{spanning tree} in a graph $G$ is a subgraph that is a tree and contains all the vertices of $G.$ 

Given $u, v\in V(G),$ we write $a_{uv}$ as the number of edges between vertices $u$ and $v.$ The matrix $A=A(G)=\{a_{uv}\}_{u, v\in V(G)}$ is called \textit{the adjacency matrix} of the graph $G.$ The degree $d_v$ of a vertex $v \in V(G)$ is defined by $d_v=\sum_{u\in V(G)}a_{uv}.$ Let $D=D(G)$ be the diagonal matrix indexed by the elements of $V(G)$ with $d_{vv} = d_v.$ The matrix $L=L(G)=D(G)-A(G)$ is called \textit{the Laplacian matrix}, or simply \textit{Laplacian}, of the graph $G.$ 

In what follows, we denote by $I_n$ the identity matrix of order $n.$

We say that an $n\times n$ matrix is \textit{circulant,} denoted by $circ(a_0, a_1,\ldots,a_{n-1}),$ if it is of the form
$$circ(a_0, a_1,\ldots, a_{n-1})=
\left(\begin{array}{ccccc}
a_0 & a_1 & a_2 & \ldots & a_{n-1} \\
a_{n-1} & a_0 & a_1 & \ldots & a_{n-2} \\
  & \vdots &   & \ddots & \vdots \\
a_1 & a_2 & a_3 & \ldots & a_0\\
\end{array}\right).$$

Recall \cite{PJDav} that the eigenvalues of the matrix $C=circ(a_0,a_1,\ldots,a_{n-1})$ are given by the following simple formulas $\lambda_j=p(\varepsilon^j_n),\,j=0,1,\ldots,n-1$ where $p(x)=a_0+a_1 x+\ldots+a_{n-1}x^{n-1}$ and $\varepsilon_n$ is an order $n$ primitive root of the unity.  For any $i=0,\ldots, n-1$, let $w_i =(1,\varepsilon_n^i ,\varepsilon_n^{2i},\ldots,\varepsilon_n^{(n-1)i})^t$ be a column vector of length $n.$ Then all $n\times n$ circulant matrices share the same set of linearly independent eigenvectors $w_0, w_1, \ldots, w_{n-1}.$ Hence, any set of $n\times n$ circulant matrices can be simultaneously diagonalizable. Moreover, let $T_{n}=circ(0,1,0,\ldots,0)$ be the matrix representation of the shift operator $T_{n}:(x_0,x_1,\ldots,x_{n-2},x_{n-1})\rightarrow(x_1, x_2,\ldots,x_{n-1},x_0).$  Then $T_{n} w_i=\varepsilon_n w_i,\,T_{n}^{-1} w_i=\varepsilon_n^{-1} w_i,$ and for any  Laurent polynomial $P(z)$ one has $P(T_{n}) w_i=P(\varepsilon_n) w_i.$

Let $D$ be a group, and let $S$ be a subset of $D,$ which doesn't contain the identity element $1.$ The Cayley digraph associated with $(D,S)$ is then defined as the directed graph with the set of vertices $D$ and the set of edges 
$$\{(g,h): g,h\in D, gh^{-1}\in S\}.$$
The Cayley graph depends on the choice of a generating set $S,$ and is connected if and only if $S$ generates $D$ (i.e., the set $S$ are group generators of $D$). We deal with undirected graphs, and always assume that $S=S^{-1},$ where $S^{-1}=\{s:s^{-1}\in S\}.$
 
Let $\mathbb{D}_n=\langle a,b|a^2=1,b^n=1,(a\,b)^2=1\rangle$ be dihedral group of order $2n.$ We arrange the elements of the group $\mathbb{D}_{n}$ as $V=\{1,b,\ldots,b^{n-1}, a, b a,\ldots,b^{n-1} a\}$ and consider the Cayley graph
\begin{equation}\label{eq:d1}\mathcal{D}_{n}=Cay(\mathbb{D}_{n}, b^{\pm\beta_1},b^{\pm\beta_2},\ldots,b^{\pm\beta_s}, a b^{\gamma_1}, a b^{\gamma_2},\ldots, a b^{\gamma_t} )\end{equation} with the generating set
$S=\{b^{\pm\beta_1},b^{\pm\beta_2},\ldots,b^{\pm\beta_s}, a b^{\gamma_1}, a b^{\gamma_2},\ldots, a b^{\gamma_t}\}$
for some integers $\beta_1,\beta_2,\cdots,$ and $\gamma_1,\gamma_2,\cdots.$ We suppose that $\mathbb{D}_{n}$ acts on $V$  by the  rule:  $g\in \mathcal{D}_{n} $ sends a vertex $v\in V$  to the vertex $v g.$ Then the set of oriented edges of $\mathcal{D}_{n}$ can be describe as follows. Given  $j\in \{\pm \beta_1,\pm \beta_2,\ldots,\pm \beta_s\}$ there is an edge  $b^{k}\xrightarrow{b^{j}}b^{k+j}$  and an edge $b^{k}a\xrightarrow{b^{j}}b^{k-j}a$  for any $k=0,1,\ldots, n-1;$  given  $j\in\{\gamma_1,\gamma_2,\ldots,\gamma_t\}$ there is an edge  $b^{k}\xrightarrow{a b^{j}}b^{k-j}a$  and an edge $b^{k}a\xrightarrow{a b^{j}}b^{k+j}$ for any $k=0,1,\ldots, n-1.$ Noting that $(a b^{j})^{-1}=b^{-j} a=a b^{j},$ we have $S=S^{-1}.$ Hence  $\mathcal{D}_{n}$ is an undirected graph.

We always restrict to the case $0<\beta_1<\beta_2<\cdots<\beta_s<\frac{n}{2}$ and $0\le\gamma_1<\gamma_2<\cdots<\gamma_t\le{n-1}.$ We suppose that $s\geq0$ and $t\geq1.$ Then the graph $\mathcal{D}_{n}$ has no loops and multiple edges.

Now we introduce a necessary and sufficient condition for the connectedness of the graph $\mathcal{D}_{n}.$

\begin{lemma}\label{gconnect}
The graph $\mathcal{D}_{n}=Cay(\mathbb{D}_{n}, b^{\pm\beta_1},b^{\pm\beta_2},\ldots,b^{\pm\beta_s}, a b^{\gamma_1}, a b^{\gamma_2},\ldots, a b^{\gamma_t} )$ is connected if and only if $\gcd(n,\,\beta_j,\,1\le j\le s,\,\gamma_j-\gamma_k,\,1\le j<k\le t)$ is equal to $1.$
\end{lemma}

\begin{pf}
Since our graph $\mathcal{D}_{n}$ is a Cayley graph associated with $(\mathbb{D}_{n},\,S),$ it is connected if and only if the set $S$ generates group $\mathbb{D}_n.$   

Suppose that $\gcd(n,\,\beta_j,\,1\le j\le s,\,\gamma_j-\gamma_k,\,1\le j<k\le t)=1.$ Then, there exist integer numbers 
$\{c_{j}\}_{1\le j\le s},\,\{d_{j,k}\}_{1\le j<k\le t},\,e$ such that 
$$\sum_{j=1}^{s}c_{j}\beta_{j}+\sum_{1\le j<k\le t}d_{j,k}(\gamma_{j}-\gamma_{k})+e\,n=1.$$ 
(See, for example, \cite{Apost}, p. 21). One easily sees that 
$$b=b^{\sum_{j=1}^{s}c_{j}\beta_{j}+\sum_{1\le j<k\le t}d_{j,k}(\gamma_{j}-\gamma_{k})+e\,n}=b^{\sum_{j=1}^{s}c_{j}\beta_{j}}b^{\sum_{1\le j<k\le t}d_{j,k}(\gamma_{j}-\gamma_{k})}$$
We rewrite $b^{\sum_{1\le j<k\le t}d_{j,k}(\gamma_{j}-\gamma_{k})}$ in the form
$$b^{\sum_{1\le j<k\le t}d_{j,k}(\gamma_{j}-\gamma_{k})}=\prod_{1\le j<k\le t}\underbrace{b^{-\gamma_{k}}b^{\gamma_{j}}\cdot\ldots\cdot b^{-\gamma_{k}}b^{\gamma_{j}}}_{d_{j,k}\text{times}}=$$
$$\prod_{1\le j<k\le t}\underbrace{b^{-\gamma_{k}}\,a\,a\,b^{\gamma_{j}}\cdot\ldots\cdot b^{-\gamma_{k}}\,a\,a\,b^{\gamma_{j}}}_{d_{j,k}\text{times}}=\prod_{1\le j<k\le t}(a\,b^{\gamma_{k}}\cdot a\,b^{\gamma_{j}})^{d_{j,k}}.$$
In other words, we can express the generator $b$ in terms of elements of $S$
$$b=\prod_{j=1}^{s}(b^{\beta_{j}})^{c_{j}}\prod_{1\le j<k\le t}(a\,b^{\gamma_{k}}\cdot a\,b^{\gamma_{j}})^{d_{j,k}}.$$
The set $S$ contains a term $a\,b^{\gamma_{1}}.$  So that $a=a\,b^{\gamma_{1}}\cdot b^{-\gamma_{1}}$ also is the product of elements of $S.$ 

On the other hand, let $\mathcal{D}_n$ be connected. Then the set $S$ generates full group $\mathbb{D}_n.$ Hence, there exist words in $S$ representing $b.$ We can always choose the representation of $b$ in the form 
$$b=\prod_{j=1}^{s}(b^{\beta_{j}})^{c_{j}}\prod_{k=1}^{N}a\,b^{\gamma_{i_k}},$$
where $N\geq0,\ c_j$ are integers for all $j=1,\ldots,s$ and each $i_k$ lies in the set $\{1,2,\ldots,t\}$ for all $k=1,\ldots,N.$ In fact, using the fact that any element $a\,b^{\gamma_j}$ is an involution and the following equalities in $\mathbb{D}_n$
$$a\,b^{\gamma_{i}}\cdot b^{\beta_{j}}=a\,b^{\beta_{j}}\cdot b^{\gamma_{i}}=b^{-\beta_{j}}\cdot a\,b^{\gamma_{i}},$$
one can reduce any word in terms of $S$ to a given form.

We note that $N$ is even. Otherwise, if $N$ is odd then using above arguments we deduce that $a$ is a proper power of $b.$ This is impossible in the dihedral group. So 
$$b=\prod_{j=1}^{s}(b^{\beta_{j}})^{c_{j}}\prod_{k=1}^{N}a\,b^{\gamma_{i_k}}=
\prod_{j=1}^{s}(b^{\beta_{j}})^{c_{j}}\prod_{k=1}^{N/2}a\,b^{\gamma_{i_{2k-1}}}\cdot a\,b^{\gamma_{i_{2k}}}=$$
$$\prod_{j=1}^{s}(b^{\beta_{j}})^{c_{j}}\prod_{k=1}^{N/2}b^{-\gamma_{i_{2k-1}}}\,a\cdot a\,b^{\gamma_{i_{2k}}}=
\prod_{j=1}^{s}(b^{\beta_{j}})^{c_{j}}\prod_{k=1}^{N/2}b^{-\gamma_{i_{2k-1}}+\gamma_{i_{2k}}}.$$

So that
$$b=b^{\sum_{j=1}^{s}{\beta_{j}}{c_{j}}+\sum_{k=1}^{N/2}(\gamma_{i_{2k}}-\gamma_{i_{2k-1}})}.$$ 
Or, (since $b$ is an element of order $n$)
$$\sum_{j=1}^{s}c_{j}\beta_{j}+\sum_{k=1}^{N/2}(\gamma_{i_{2k}}-\gamma_{i_{2k-1}})\equiv1(\mod n).$$
Hence, there exists an integer $e$ such that 
$$\sum_{j=1}^{s}c_{j}\beta_{j}+\sum_{k=1}^{N/2}(\gamma_{i_{2k}}-\gamma_{i_{2k-1}})+e\,n=1.$$
Therefore,  $\gcd(n,\,\beta_j,\,1\le j\le s,\,\gamma_j-\gamma_k,\,1\le j<k\le t)=1$
and the lemma is proved.\end{pf}

The  adjacency matrix of the graph $\mathcal{D}_{n}$ is given by the $2n\times2n$  block matrix
$$A=\left(\begin{array}{cc}\sum\limits_{j=1}^s(T_n^{\beta_j}+T_n^{-\beta_j})&\sum\limits_{j=1}^tT_n^{-\gamma_j}\\ \sum\limits_{j=1}^tT_n^{\gamma_j}&\sum\limits_{j=1}^s(T_n^{\beta_j}+T_n^{-\beta_j})\end{array}\right),$$ where $T_n=circ (0,1,0,\ldots,0).$  The corresponding degree matrix is $D=\left(\begin{array}{cc}(2s+t)I_n&0\\ 0&(2s+t)I_n\end{array}\right),$ where $I_n$ is the $n\times n$ identity matrix.  Since the Laplacian of $\mathcal{D}_n$ is  $L=D-A,$  we  have
$$L=\left(\begin{array}{cc}(2s+t)I_n-\sum\limits_{j=1}^s(T_n^{\beta_j}+T_n^{-\beta_j})&-\sum\limits_{j=1}^tT_n^{-\gamma_j}\\ -\sum\limits_{j=1}^tT_{n}^{\gamma_j}&(2s+t)I_n-\sum\limits_{j=1}^s(T_n^{\beta_j}+T_n^{-\beta_j})\end{array}\right).$$ 

\section{Spectrum of a Cayley graph on the dihedral group}\label{spectra}

Let $\mathcal{D}_{n}=Cay(\mathbb{D}_{n}, b^{\pm\beta_1},b^{\pm\beta_2},\ldots,b^{\pm\beta_s}, a b^{\gamma_1}, a b^{\gamma_2},\ldots, a b^{\gamma_t} )$ be the Cayley graph on a dihedral group $\mathbb{D}_{n}.$ We introduce the following Laurent polynomials 
$${\cal A}(z)=2s+t-\sum\limits_{i=1}^{s}(z^{\beta_{i}}+z^{-\beta_{i}}),\,{\cal B}(z)=-\sum\limits_{i=1}^{t}z^{\gamma_{i}}
\text{ and }P(z)={\cal A}(z){\cal A}(z^{-1})-{\cal B}(z){\cal B}(z^{-1}).$$ We will refer to $P(z)$ as \textit{the Laurent polynomial associated with the graph} $\mathcal{D}_{n}.$ Then the Laplacian of the graph $\mathcal{D}_{n}$ is given by the following $2n\times2n$ block matrix $L=\left(\begin{array}{cc}{\cal A}(T)&{\cal B}(T^{-1})\\{\cal B}(T)&{\cal A}(T^{-1})\end{array}\right),$ where $T=circ(0,1,0,\ldots,0).$ 

Recall that $T$ is conjugated to $\mathbb{T}=diag(1,\varepsilon,\ldots,\varepsilon^{n-1}),$ where $\varepsilon=\varepsilon_{n}=\exp(2\pi\mathtt{i}/n).$ So that $L$ is conjugated to the matrix $\mathbb{L}=\left(\begin{array}{cc}{\cal A}(\mathbb{T})&{\cal B}(\mathbb{T}^{-1})\\
{\cal B}(\mathbb{T})&{\cal A}(\mathbb{T}^{-1})\end{array}\right).$ The spectra of $L$ and $\mathbb{L}$ are the same. It can be found by solving the system of linear equations 
$$\begin{cases}
{\cal A}(\mathbb{T})\textbf{x}+{\cal B}(\mathbb{T}^{-1})\textbf{y}=\lambda \textbf{x}\\ 
{\cal B}(\mathbb{T})\textbf{x}+{\cal A}(\mathbb{T}^{-1})\textbf{y}=\lambda \textbf{y},
\end{cases}$$
where $\textbf{x},\textbf{y}\in\mathbb{C}^{n}$ and $(\textbf{x},\textbf{y})\neq(\textbf{0},\textbf{0}).$ Since the matrix $\mathbb{T}$ is diagonal, the system of equations splits into $n$ scalar linear systems 
$$\begin{cases}
{\cal A}(\varepsilon^{j})x+{\cal B}(\varepsilon^{-j})y=\lambda x\\
{\cal B}(\varepsilon^{j})x+{\cal A}(\varepsilon^{-j})y=\lambda y,
\end{cases}$$ 
where $j=0,1\ldots,n-1.$ Hence, $\lambda$ is a root of the quadratic equation  
$$\lambda^2-({\cal A}(\varepsilon^{j})+{\cal A}(\varepsilon^{-j}))\lambda+
{\cal A}(\varepsilon^{j}){\cal A}(\varepsilon^{-j})-{\cal B}(\varepsilon^{j}){\cal B}(\varepsilon^{-j})=0.$$
The solutions of this equation are $\lambda_{j,1}=\Re({\cal A}(\varepsilon^{j}))+\sqrt{-\Im({\cal A}(\varepsilon^{j}))^{2}+
|{\cal B}(\varepsilon^{j})|^{2}}$ and $\lambda_{j,2}=\Re({\cal A}(\varepsilon^{j}))-
\sqrt{-\Im({\cal A}(\varepsilon^{j}))^{2}+|{\cal B}(\varepsilon^{j})|^{2}}.$ We note that L.~Babai  \cite{Babai} found similar formulas for eigenvalues for the adjacency matrices of $\mathcal{D}_n$ by making use of the representation theory for finite groups. The corresponding eigenvectors of the operator $\mathbb{L}$ are $u_{j,1}={\cal B}(\varepsilon^{j})\textbf{e}_{j+1}+(\lambda_{j,1}
-{\cal A}(\varepsilon^{j}))\textbf{e}_{j+1}$ and $u_{j,2}={\cal B}(\varepsilon^{j})\textbf{e}_{j+1}+(\lambda_{j,2}
-{\cal A}(\varepsilon^{j}))\textbf{e}_{j+1},$ where $\textbf{e}_{j}$ is the $j$-th basic vector in the $\mathbb{C}^{2n}.$ 

For the special case $j=0,$ we have $\lambda^2-2{\cal A}(1)\lambda+{\cal A}(1)^2-{\cal B}(1)^2=0.$ Hence $$\lambda_{0,1}={\cal A}(1)+{\cal B}(1)=0,\,\lambda_{0,2}={\cal A}(1)-{\cal B}(1)=2{\cal A}(1).$$ We note also that $\lambda_{j,1}\lambda_{j,2}={\cal A}(\varepsilon^{j}){\cal A}(\varepsilon^{-j})-{\cal B}(\varepsilon^{j}){\cal B}(\varepsilon^{-j})=P(\varepsilon^{j}).$ We suppose that the graph $\mathcal{D}_{n}$ is connected. Then all the Laplacian eigenvalues, except of $\lambda_{0,1}=0,$ are non-zero. They are 
$$\lambda_{0,2},\,\lambda_{j,1},\,\lambda_{j,2},\,j=1,2,\ldots,n-1.$$ By the Kirchhoff's theorem we have the following formula for the number of spanning trees in $\mathcal{D}_{n}.$ 

$$\tau(n)=\frac{\lambda_{0,2}}{2n}\prod\limits_{j=1}^{n-1}\lambda_{j,1}\lambda_{j,2}
=\frac{{\cal A}(1)}{n}\prod\limits_{j=1}^{n-1}P(\varepsilon^{j}).$$ Since ${\cal A}(1)=t,$ we get the following result.  

\begin{theorem}\label{prelimenary}
Let $\mathcal{D}_{n}=Cay(\mathbb{D}_{n}, b^{\pm\beta_1},b^{\pm\beta_2},\ldots,b^{\pm\beta_s}, a b^{\gamma_1}, a b^{\gamma_2},\ldots, a b^{\gamma_t} )$ be a Cayley graph on the dihedral group $\mathbb{D}_{n}.$ Then the number of spanning trees of $\mathcal{D}_{n}$ is given by the formula $$\tau(n)=\frac{t}{n}\prod\limits_{j=1}^{n-1}P(\varepsilon_{n}^{j}),$$ where $\varepsilon_{n}=\exp(2\pi\mathtt{i}/n),\,P(z)={\cal A}(z){\cal A}(z^{-1})
-{\cal B}(z){\cal B}(z^{-1}),\,{\cal A}(z)=2s+t-\sum\limits_{i=1}^{s}(z^{\beta_{i}}+z^{-\beta_{i}})$ and
${\cal B}(z)=-\sum\limits_{i=1}^{t}z^{\gamma_{i}}.$
\end{theorem}

The following two lemmas give more properties of the associated polynomial $P(z).$

\begin{lemma}\label{lemma1} Let the graph $\mathcal{D}_{n}$ be connected. Then $P(1)=0, P^{\prime}(1)=0$ and $P^{\prime\prime}(1)<0.$ \end{lemma}
\begin{pf} By direct calculation we have $P(1)=P^{\prime}(1)=0$ and  
$$P^{\prime\prime}(1)=-4t\sum_{j=1}^{s}\beta_{j}^2+2(\sum_{j=1}^{t}\gamma_{j})^2-2t\sum_{j=1}^{t}\gamma_{j}^2.$$
We use Lagrange's identity 
$$(\sum_{j=1}^{t}\xi_{j}^2)(\sum_{j=1}^{t}\zeta_{j}^2)-(\sum_{j=1}^{t}\xi_{j}\zeta_{j})^2=\sum\limits_{1\le j<k\le t}(\xi_{j}\zeta_{k}-\xi_{k}\zeta_{j})^2.$$ 
Setting $\xi_j=1$ and $\zeta_j=\gamma_j,$ we get 
$$P^{\prime\prime}(1)=-2(2t\sum_{j=1}^{s}\beta_{j}^2+\sum_{1\le j<k\le t}(\gamma_j-\gamma_k)^2)<0.$$
In fact, it is impossible that all elements $\{\beta_{i}\}_{i=1,\ldots,s}$ and $\{\gamma_{j}-\gamma_{k}\}_{j,k=1,\ldots,t}$ are equal to $0,$ since the graph $\mathcal{D}_{n}$ is connected. See Lemma~\ref{gconnect}. 
\end{pf}  
\begin{lemma}\label{lemma2} Suppose that the numbers $\{\beta_j,\,1\le j\le s,\,\gamma_j-\gamma_k,\,1\le j<k\le t\}$  are relatively   prime. Then for any  $\varphi\in\mathbb{R}$, we have $P(e^{\mathtt{i} \varphi})\ge0.$  Furthermore,  $P(e^{\mathtt{i} \varphi})=0$ if and
only if $e^{\mathtt{i} \varphi}=1.$ 
 \end{lemma}
 \begin{pf}  We have $P(e^{\mathtt{i} \varphi})={\cal A}(e^{\mathtt{i} \varphi}){\cal A}(e^{-\mathtt{i} \varphi})-{\cal B}(e^{\mathtt{i} \varphi}){\cal B}(e^{-\mathtt{i} \varphi})=|{\cal A}(e^{\mathtt{i} \varphi})|^2-|{\cal B}(e^{\mathtt{i} \varphi})|^2.$ Also,
$$
|A(e^{\mathtt{i} \varphi})|^2=(s+t-\sum\limits_{j=1}^s\cos \beta_j\varphi)^2+(\sum\limits_{j=1}^s\sin \beta_j\varphi)^2
= (t+\sum\limits_{j=1}^s(1-\cos \beta_j\varphi))^2+ (\sum\limits_{j=1}^s\sin \beta_j\varphi)^2\ge t^2
$$ and
 $$|{\cal B}(e^{\mathtt{i} \varphi})|^2=|\sum_{j=1}^t e^{\mathtt{i}\gamma_j}|^2=\sum\limits_{1\le j,\,k\le t}\cos(\gamma_j-\gamma_k)\varphi=
t+2\sum\limits_{1\le j<k\le t}\cos(\gamma_j-\gamma_k)\varphi\le t^2.$$
 
Hence, $P(e^{\mathtt{i} \varphi})=|{\cal A}(e^{\mathtt{i} \varphi})|^2-|{\cal B}(e^{\mathtt{i} \varphi})|^2\ge 0.$ We have the equality if and only if 
$$(t+\sum\limits_{j=1}^s(1-\cos \beta_j\varphi))^2+(\sum\limits_{j=1}^s\sin \beta_j\varphi)^2=t^2 \text{  and  } t^2=t+2\sum\limits_{1\le j<k\le t}\cos(\gamma_j-\gamma_k)\varphi.$$  
The first equation  is equivalent to $$\cos \beta_j\varphi=1,1\le j\le s,$$ and the second one is equivalent to
$$\cos(\gamma_j-\gamma_k)\varphi=1,\,1\le j<k\le t.$$
 
Since the numbers $\{\beta_j,\,1\le j\le s,\,\gamma_j-\gamma_k,\,1\le j<k\le t\}$ are relatively prime, we obtain $\cos \varphi =1,$  or $e^{\mathtt{i} \varphi}=1.$ 
\end{pf}

\section{Counting spanning trees}\label{counting}

The main results of this section paper are the following theorems.

\begin{theorem}\label{theorem1} Let 
$\mathcal{D}_{n}=Cay(\mathbb{D}_{n}, b^{\pm\beta_1},b^{\pm\beta_2},\ldots,b^{\pm\beta_s}, a b^{\gamma_1}, a b^{\gamma_2},\ldots, a b^{\gamma_t} )$ be a Cayley graph on the  group $\mathbb{D}_{n}$ and $P(z)$ be the associated Laurent polynomial for $\mathbb{D}_{n}.$ Then the number of spanning trees $\tau(n)$ in the graph $\mathcal{D}_{n}$ is given by the formula $$\tau(n)= \frac{n\,t\,|\eta|^{n}}{q}\prod\limits_{\substack{P(z)=0\\z\neq1}}|z^n-1|,$$
where the product is taken over all the roots, different from $1,$ of the associated Laurent polynomial $P(z),$ {} $\eta$ is the leading coefficient of $P(z),\text{ and } q=2t\sum_{j=1}^{s}\beta_{j}^2+\sum_{1\le j<k\le t}(\gamma_j-\gamma_k)^2.$ \end{theorem}

\begin{pf}
By Theorem \ref{prelimenary} we already have
\begin{equation}\label{tauHnew}
\tau(n)=\frac{t}{n}\prod\limits_{j=1}^{n-1}P(\varepsilon_n^j).
\end{equation}
Denote by $\eta$ the leading coefficient of $P(z).$ Since $P(z)=P(1/z)$ we can present the polynomial in the form 
$$P(z)=\eta\,z^{-r}+a_{1}z^{-r+1}+\ldots+a_{r}+\ldots+a_{1}z^{r-1}+\eta\,z^{r},$$   
for some $r\geq0.$ To continue the proof, we replace the Laurent polynomial $P(z)$ by $\widetilde{P}(z)=\frac{z^{r}}{\eta}P(z).$ Then $\widetilde{P}(z)$ is a monic polynomial of the degree $2r$ with the same roots as $P(z).$ We note that
\begin{equation}\label{newPnew}
\prod\limits_{j=1}^{n-1}\widetilde{P}(\varepsilon_{n}^{j})=\frac{\varepsilon_{n}^{\frac{(n-1)n}{2}r}}{\eta^{n-1}} \prod\limits_{j=1}^{n-1}P(\varepsilon_{n}^{j})=\frac{(-1)^{r(n-1)}}{\eta^{n-1}}\prod\limits_{j=1}^{n-1}P(\varepsilon_{n}^{j}).
\end{equation}

By Lemma \ref{lemma1} the polynomial $\widetilde{P}(z)$ has two roots equal to $1$ and all the other roots different from $1.$ Also, we recognize the complex numbers $\varepsilon_{n}^{j},\,j=1,\ldots,n-1$ as the roots of polynomial $\frac{z^n-1}{z-1}.$ By the basic properties of resultant we have (\cite{Pras}, Ch. 1.3)
\begin{eqnarray}\label{Hlemmanew}
\nonumber &&\prod\limits_{j=1}^{n-1}\widetilde{P}(\varepsilon_{n}^{j})=\textrm{Res}(\widetilde{P}(z),\frac{z^{n}-1}{z-1}) =\textrm{Res}(\frac{z^{n}-1}{z-1},\widetilde{P}(z))=\prod\limits_{z:\,\widetilde{P}(z)=0}\frac{z^{n}-1}{z-1}\\
&&=\prod\limits_{z:\,P(z)=0}\frac{z^{n}-1}{z-1}=(\lim_{z\to1}\frac{z^n-1}{z-1})^{2}\prod\limits_{\substack{P(z)=0\\z\neq1}}\frac{z^n-1}{z-1}= n^{2}\prod\limits_{\substack{P(z)=0\\z\neq1}}\frac{z^n-1}{z-1}.
\end{eqnarray}
Combining (\ref{tauHnew}), (\ref{newPnew}) with (\ref{Hlemmanew}), we have the following formula for the number of spanning trees 
\begin{eqnarray}\label{before}\tau(n)=(-1)^{r(n-1)}\eta^{n-1}n\,t\prod\limits_{\substack{P(z)=0\\z\neq1}}\frac{z^n-1}{z-1}.
\end{eqnarray} 

To finish the proof we need to evaluate the product $\prod\limits_{\substack{P(z)=0\\z\neq1}}(z-1)=\prod\limits_{\substack{
\widetilde{P}(z)=0\\z\neq1}}(z-1).$ To do this, we use the following property of resultants
(\cite{Pras}, Ch. 1.3, property 3). If $B$ is a monic polynomial, then for any polynomials $A$ and $C$ we have 
$$\textrm{Res}(B,A+CB)=\textrm{Res}(B,A).$$
We use the Taylor expansion of polynomial $\widetilde{P}(z)$ at the point $z=1,$ which has the form 
$$\widetilde{P}(z)=\widetilde{P}(1)+\widetilde{P}^{\prime}(1)(z-1)+\frac{\widetilde{P}^{\prime\prime}(1)}{2}(z-1)^{2}+\frac{\widetilde{P}^{\prime\prime\prime}(1)}{6}(z-1)^{3}+\ldots.$$
Since $\widetilde{P}(1)=\widetilde{P}^{\prime}(1)=0,$ we have 
$$\prod\limits_{\substack{\widetilde{P}(z)=0\\z\neq1}}(z-1)=\textrm{Res}(z-1,
\frac{\widetilde{P}(z)}{(z-1)^2})=\textrm{Res}(z-1,\frac{\widetilde{P}^{\prime\prime}(1)}{2}+\frac{\widetilde{P}^{\prime\prime\prime}(1)}{6}(z-1)+\ldots)=\textrm{Res}(z-1,\frac{\widetilde{P}^{\prime\prime}(1)}{2}).$$ The last term is equal to 
$\frac{\widetilde{P}^{\prime\prime}(1)}{2}=\frac{1}{2}(\frac{z^{r}}{\eta}P(z))^{\prime\prime}_{z=1}=\frac{P^{\prime\prime}(1)}{2 \eta}=-\frac{q}{\eta}$ (see Lemma~\ref{lemma1}). As a result, we get
\begin{eqnarray}\label{after}\tau(n)=\frac{(-1)^{r(n-1)+1}n\,t\,\eta^{n}}{q}\prod\limits_{\substack{P(z)=0\\z\neq1}}(z^n-1).\end{eqnarray} Since $\tau(n)$ is a positive integer, the statement of theorem follows.
\end{pf} 

\noindent\textbf{Remark.} According to Lemma \ref{lemma1}, for the calculation of $q$ one can get the formula $q=-\frac{1}{2}P^{\prime\prime}(1).$
\medskip

The associated polynomial $P(z)$ is a palindromic Laurent polynomial, satisfying the property $P(z)=P(1/z).$ In other words, $P(z)$ has the form $P(z)=c_{0}+\sum_{k=1}^{r}c_{k}(z^{k}+z^{-k}).$  Denote by $T(n,x)=\cos(n\arccos(x))$ the $n$-th Chebyshev polynomial of the first kind. The following equality is known $T(n,\frac{z+z^{-1}}{2})=\frac{z^{n}+z^{-n}}{2}.$ Other basic properties of Chebyshev polynomials can be found in \cite{MasonHandCheb}. 

From here, we conclude that $$P(z)=Q(\frac{z+z^{-1}}{2}),\text{ where }Q(w)=c_{0}+\sum_{k=1}^{r}2c_{k}T(k,w).$$  
We refer to the polynomial $Q(w)$ as \textit{the Chebyshev transform} of $P(z).$ These two polynomials are also related as $Q(w)=P(w+\sqrt{w^{2}-1}).$ It is easier to deal with $Q(w)$ since it is an ordinary polynomial with degree twice less than $P(z).$ 

The following theorem is a direct consequence of Theorem~\ref{theorem1}. 

\begin{theorem}\label{theorem2} Let $P(z)$ be the associated Laurent polynomial of the Cayley graph $\mathcal{D}_n.$ Denote by $Q(w)$  the Chebyshev transform of $P(z)$ and let $r$ be the degree of polynomial $Q(w).$  The number of spanning trees $\tau(n)$ in the graph $\mathcal{D}_n$ is given by the formula 
$$\tau(n)=\frac{n\,t |\eta|^n}{q}\prod_{p=1}^{r-1}|2T(n,w_p)-2|,$$
where $w_p,\,p=1,2,\ldots,r-1$ are roots, different from $1,$ of the algebraic equation $Q(w)=0,$ and $T(n,w)$ is the Chebyshev polynomial of the first kind, $\eta$ is the leading coefficient of $P(z),$ and $q=2t\sum_{j=1}^{s}\beta_{j}^2+\sum_{1\le j<k\le t}(\gamma_j-\gamma_k)^2.$
\end{theorem}

\begin{pf} Since $P(z)=Q((z+1/z)/2),$  by Lemma~\ref{lemma1}, the roots of polynomials ${P}(z)$ and $Q(w)$ are $1,1,z_{1},1/z_{1},\ldots,z_{r-1},1/z_{r-1},\,z_{j}\neq1\textrm{ and }1,\,\frac{1}{2}(z_{j}+z_{j}^{-1})=w_{j}\neq1,\,j=1,\ldots,r-1,$ respectively. Also  $T(n, w_{j})=\frac{z_{j}^{n}+z_{j}^{-n}}{2}.$ 
Hence, $$\prod\limits_{\substack{P(z)=0\\z\neq1}}(z^n-1)=\prod\limits_{j=1}^{r-1}({z_{j}^{n}-1})({z_{j}^{-n}-1})= \prod\limits_{j=1}^{r-1}(2-{z_{j}^{n}-z_{j}^{-n}})=-\prod\limits_{j=1}^{r-1} ({2T(n, w_{j})-2}).$$ Substituting  the latter in Theorem~\ref{theorem1}
we get the result.
\end{pf}
\section{Arithmetical properties of complexity for  the graph $\mathcal{D}_n$}\label{arithmetic}

The main result of this section is the following theorem.
 
\bigskip

\begin{theorem}\label{lorenzini} Let $\tau(n)$ be the number of spanning trees in the graph 
$$\mathcal{D}_{n}=Cay(\mathbb{D}_{n}, b^{\pm\beta_1},b^{\pm\beta_2},\ldots,b^{\pm\beta_s}, a b^{\gamma_1}, a b^{\gamma_2},\ldots, a b^{\gamma_t} ).$$
Denote by $\beta_{odd}$ ($\gamma_{odd}$ resp.) the quantity of odd numbers in the sequences $\beta_{1}, \beta_{2},\ldots, \beta_{s}$ ($\gamma_{1}, \gamma_{2}, \ldots, \gamma_{t}$ resp.). Also, denote by $\gamma_{even}$ the quantity of even numbers in the sequence $\gamma_{1}, \gamma_{2}, \ldots, \gamma_{t}.$ Let $\delta$ be the square free part of the integer $\xi=(2\beta_{odd}+\gamma_{odd})(2\beta_{odd}+\gamma_{even}).$ Then there exists an integer sequence $a(n)$ such that
\begin{enumerate}
\item[$1^\circ$] $\tau(n)= n\,t\,a(n)^{2},$ if $n$ is odd;
\item[$2^\circ$] $\tau(n)= n\,t\,\delta\,a(n)^{2},$ if $n$ is even.
\end{enumerate}

\end{theorem}

\begin{pf} We consider the associated polynomial $P(z)$ for the graph $\mathcal{D}_{n}.$ Now we will express the value $P(-1)$ through basic parameters of $\mathcal{D}_{n}.$ Denote by $\beta_{even}$ ($\beta_{odd}$ resp.) the number of even numbers (odd numbers resp.) in the sequence $\beta_{1}, \beta_{2},\ldots, \beta_{s}.$ Denote by $\gamma_{even}$ ($\gamma_{odd}$ resp.) the number of even numbers (odd numbers resp.) in the sequence $\gamma_{1}, \gamma_{2}, \ldots, \gamma_{t}.$ One easily sees that $s=\beta_{even}+\beta_{odd}$ and $t=\gamma_{even}+\gamma_{odd}.$ Since $P(z)=\mathcal{A}(z)\mathcal{A}(z^{-1})-\mathcal{B}(z)\mathcal{B}(z^{-1}),$ we get 
\begin{eqnarray*}
P(-1)&=&\mathcal{A}(-1)^2-\mathcal{B}(-1)^2=(2s+t-\sum_{j=1}^{s}((-1)^{\beta_{j}}+(-1)^{-\beta_{j}}))^{2}-(\sum_{j=1}^{t}(-1)^{-\gamma_{j}})^{2}\\
&=&(2s+t-2\sum_{j=1}^{s}(-1)^{\beta_{j}})^{2}-(\sum_{j=1}^{t}(-1)^{-\gamma_{j}})^{2}\\
&=&(2(\beta_{even}+\beta_{odd})+\gamma_{even}+\gamma_{odd}-2(\beta_{even}-\beta_{odd}))^2-(\gamma_{even}-\gamma_{odd})^2\\
&=&4(2\beta_{odd}+\gamma_{even})(2\beta_{odd}+\gamma_{odd}).
\end{eqnarray*}
As a consequence, we have $P(-1)=4\xi=\delta(2\omega)^2,$ where $\delta$ is the square free part of $\xi$ and $\omega$ is some integer.

By formula~(\ref{tauHnew}) we have $n\,\tau(n)=t\prod_{j=1}^{n-1}\lambda_{j,1}\lambda_{j,2}.$ Note that $\lambda_{j,1}\lambda_{j,2}=P(\varepsilon_{n}^{j})= P(\varepsilon_{n}^{n-j})=\lambda_{n-j,1}\lambda_{n-j,2}.$ Define $c(n)=\prod\limits_{j=1}^{\frac{n-1}{2}}\lambda_{j,1}\lambda_{j,2},$ if $n$ is odd, and $d(n)=\prod\limits_{j=1}^{\frac{n}{2}-1}\lambda_{j,1}\lambda_{j,2},$ if $n$ is even. Following \cite{Lor} we note that each algebraic number $\lambda_{i,j}$ comes into both products $\prod_{j=1}^{(n-1)/2}\lambda_{j,1}\lambda_{j,2}$ and $\prod_{j=1}^{n/2-1}\lambda_{j,1}\lambda_{j,2}$ with all of its Galois conjugate elements. Therefore, both products $c(n)$ and $d(n)$ are integer numbers. Moreover, if $n$ is even, then we get $\lambda_{\frac{n}{2},1}\lambda_{\frac{n}{2},2}=P(-1).$ Now we have
\begin{itemize}
\item[$1^\circ$] $n\,\tau(n)=t\,c(n)^2$ if $n$ is odd, 
\item[$2^\circ$] $n\,\tau(n)=t\,P(-1)\,d(n)^2=4t\,\delta\,\omega^2\,d(n)^2$ if $n$ is even.   
\end{itemize}

Note that using formula (\ref{before}) from the proof of Theorem \ref{theorem1} we conclude that $\frac{\tau(n)}{n\,t}$ is an integer. Indeed, since $\widetilde{P}(z)$ and $P(z)$ share the roots we have
$$\tau(n)=(-1)^{r(n-1)}\eta^{n-1}n\,t\prod\limits_{\substack{\widetilde{P}(z)=0\\z\neq1}}\frac{z^n-1}{z-1}.$$ 
The last product is equal to resultant of two integer polynomials $\frac{\widetilde{P}(z)}{(z-1)^2}$ and $\frac{z^n-1}{z-1}$ and, hence, it is an integer number. So that $\frac{\tau(n)}{n\,t}$ is also an integer. We get
\begin{itemize}
\item[$1^\circ$] $\frac{\tau(n)}{n\,t}=(\frac{c(n)}{n})^2$ if $n$ is odd, 
\item[$2^\circ$] $\frac{\tau(n)}{n\,t}=\delta\,(\frac{2\omega\,d(n)}{n})^2$ if $n$ is even.
\end{itemize}

As $\frac{\tau(n)}{n\,t}$ is an integer and $\delta$ is square free, all squared rational numbers in $1^{\circ}$ and $2^{\circ}$ are integers. We set $a(n)=\frac{c(n)}{n}$ if $n$ is odd and $a(n)=\frac{2\omega\,d(n)}{n}$ if $n$ is even. This proves the theorem. \end{pf}

\section{Asymptotic formulas for the number of spanning trees}\label{asymptotic}

In this section we find the asymptotics for the number of spanning trees in the graph 
$$\mathcal{D}_{n}=Cay(\mathbb{D}_{n}, b^{\pm\beta_1},b^{\pm\beta_2},\ldots,b^{\pm\beta_s}, a b^{\gamma_1}, a b^{\gamma_2},\ldots, a b^{\gamma_t} ).$$
To do this we suppose that parameters $\beta_{1},\ldots,\beta_{s}, \gamma_{1},\ldots,\gamma_{t}$ are fixed, and the inequalities $0<\beta_1<\beta_2<\cdots<\beta_s<\frac{n}{2}$ and $0\le\gamma_1<\gamma_2<\cdots<\gamma_t\le{n-1}$  hold for all sufficiently large values $n.$ We suppose also that graphs  $\mathcal{D}_{n}$ are connected.
  
According to Lemma~\ref{gconnect} we have $(n,d)=1,$ where $d=\gcd(\beta_j,\,1\le j\le s,\,\gamma_j-\gamma_k,\,1\le j<k\le t).$ Since $n$ is an arbitrarily large number and $d$ is fixed, we can choose $n$ to be a multiple of $d.$ Then $d=(n,d)=1.$ 
So that the condition of Lemma \ref{lemma2} is satisfied. This gives the following result.

\bigskip
\begin{theorem}\label{theorem3} Let $\mathcal{D}_{n}=Cay(\mathbb{D}_{n}, b^{\pm\beta_1},b^{\pm\beta_2},\ldots,b^{\pm\beta_s}, a b^{\gamma_1}, a b^{\gamma_2},\ldots, a b^{\gamma_t} )$ be an infinite family of connected graphs. Then the asymptotic behaviour for the number of spanning trees $\tau(n)$ for the graph $\mathcal{D}_{n}$ is given by the formula $$\tau(n)\sim \frac{n\,t}{q}A^n,\,n\to\infty,$$ where $A=\exp\left({\int\limits_{0}^{1}\log P(e^{2\pi\texttt{i}t})\textrm{d}t}\right)$ and $q=2t\sum_{j=1}^{s}\beta_{j}^2+\sum_{1\le j<k\le t}(\gamma_j-\gamma_k)^2.$
\end{theorem}

\begin{pf}
By Theorem \ref{theorem1} we have $\tau(n)=\frac{n\,t\,|\eta|^{n}}{q}\prod\limits_{j=1}^{r-1}|{2T(n,w_{j})-2}|,$ where  $w_{j},\,j=1,2,\ldots,r-1$ are roots, different from $1,$ of the Chebyshev transform of $P(z).$

By Lemma~\ref{lemma2}, $T(n,w_{j})=\frac{z_{j}^{n}+z_{j}^{-n}}{2},$ where the $z_{j}$ and $1/z_{j}$ are roots of the polynomial $P(z)$ with the property $|z_{j}|\neq1,\,j=1,2,\ldots,r-1.$ Replacing $z_{j}$ by $1/z_{j},$ if it is necessary, we can assume that  $|z_j|>1$ for all $j=1,2,\ldots,r-1.$ Then $T(n,w_{j})\sim\frac{1}{2}z_{j}^{n}$ and $|2T(n,w_{j})-2|\sim|z_{j}|^{n}$ as $n\to\infty.$ Hence
$$\frac{n\,t|\eta|^{n}}{q}\prod_{j=1}^{r-1}|2T_{n}(w_{j})-2|\sim\frac{n\,t|\eta|^{n}}{q}\prod_{j=1}^{r-1}|z_{j}|^{n}=\frac{n\,t}{q}|\eta|^{n}\prod\limits_{\substack{P(z)=0,\\ |z|>1}}|z|^{n}=\frac{n\,t\,A^n}{q},$$ where $A=|\eta |\prod\limits_{P(z)=0,\, |z|>1}|z|$ is the Mahler measure of the polynomial $P(z).$ By (\cite{EverWard}, p.~67), we have $A=\exp\left(\int_{0}^{1}\log|P(e^{2 \pi i t })|\textrm{d}t\right).$ The theorem is proved.
\end{pf}

\section{Generating function for the number of spanning trees}\label{genfuction}

In this section, our aim is to prove the following result.
\bigskip

\begin{theorem}\label{theoremR1} Let $\tau(n)$ be the number of spanning trees in the graph $\mathcal{D}_{n}.$ Then $F(x)=\sum\limits_{n=1}^\infty\tau(n)x^n$ is a rational function with integer coefficients. Moreover, $F(\eta\,x) = F(\frac{1}{\eta\,x}),$ where $\eta$ is the leading coefficient of the associated polynomial $P(z).$ The latter allows to represent $F(x)$ as a rational function of $u=\frac{1}{2}(\eta \,x+\frac{1} {\eta\, x} ).$
\end{theorem}
\bigskip
The proof of Theorem \ref{theoremR1} is based on the following proposition proved in \cite{MedMed2020R}.

\begin{prop}\label{proposition1} Let $R(z)$  be a degree $2s$ polynomial with integer coefficients.  Suppose that all  the roots  of the polynomial  $R(z)$  are $\xi_1,\xi_2,\ldots,\xi_{2s-1},\xi_{2s}.$ Then 
$$F(x)=\sum\limits_{n=1}^\infty\Big(n\prod\limits_{j=1}^{2s}(\xi_{j}^{n}-1)\Big)x^n$$ 
is a rational  function with integer coefficients.  

Moreover, if  $\xi_{j+s}=\xi_{j}^{-1},\,j=1,2,\ldots,s,$  then $F(x)=F(1/x).$
\end{prop}
\bigskip
Proof of Theorem \ref{theoremR1}.  By formula (\ref{after}) we have

 $$F(x)= \sum\limits_{n=1}^\infty\tau(n)x^n=  \sum\limits_{n=1}^\infty\Big(\frac{(-1)^{r(n-1)+1}n\,t\,\eta^{n}}{q}\prod\limits_{\substack{P(z)=0\\z\neq1}}(z^n-1)\Big)x^n.$$ Since all the roots of $P(z),$ different from $1,$ are
 $z_{1},1/z_{1},\ldots,z_{r-1},1/z_{r-1},$ we can rewrite the latter as

$$F(x)=\frac{(-1)^{-r+1}t}{q} \sum\limits_{n=1}^{\infty}\Big(n\prod\limits_{j=1}^{r-1}(z_j^n-1)(z_j^{-n}-1)\Big)((-1)^r \eta\,x)^n.$$ 
Since $t$ and $q$ are rational numbers, by Proposition \ref{proposition1}, $F(x)$ is a rational function with 
integer coefficients satisfying $F((-1)^r\eta\, x) = F( \frac{1 }{(-1)^r \eta\,x}).$ Hence, $F(\eta\,x) = F( \frac{1 }{\eta\,x}).$

\section{Examples}\label{example}

\subsection{Prism graph $\mathcal{D}_n=Cay(\mathbb{D}_{n},b^{\pm 1},a).$}\label{example1}

\noindent$\textbf{1}^{\circ}.$ \textbf{The number of spanning trees.} 
The associated Laurent polynomial and its Chebyshev transform are 
$$P(z)=z^{-2}-6z^{-1}+10-6z+z^2\text{ and }Q(w)=4(w-2)(w-1).$$
Here $t=1,\,\eta=1,\,q=2.$ Hence, by Theorem~\ref{theorem2} we have
$$\tau(n)=\frac{n\,t\,\eta^{n}}{q}(2T(n,2)-2)=n(T(n,2)-1).$$
This coincides with the well-known result in \cite{BoePro}.\bigskip

\noindent$\textbf{2}^{\circ}.$ \textbf{The asymptotics of $\tau(n).$}
By Theorem~\ref{theorem3}, $\tau(n)\cong \frac{n}{2} A^n,$ where $A=2+\sqrt{3}.$\bigskip

\noindent$\textbf{3}^{\circ}.$ \textbf{The generating function of $\tau(n).$}
Theorem~\ref{theoremR1} gives
$$F(x)=\sum\limits_{n=1}^\infty\tau(n)x^n=\frac{-3 + u + u^2}{2 (-2 + u)^2 (-1 + u)},$$  
where $u=\frac{1}{2}(x+\frac{1}{x}).$\bigskip

\noindent$\textbf{4}^{\circ}.$ \textbf{Divisibility by squares.}
To see the divisibility by squares, consider a few terms of generating function
$$F(x)=x + 12 x^2 + 75 x^3 + 384 x^4 + 1805 x^5 + 8100 x^6 + 35287 x^7 + 150528 x^8 +\ldots .$$
By Theorem \ref{lorenzini}, we have $\xi=(2\beta_{odd}+\gamma_{odd})(2\beta_{odd}+\gamma_{even})=(2\cdot1)(2\cdot1+1)=6.$ Hence, $\delta=6.$ So that there exists an integer sequence $a(n)$ such that $\tau(u)=n\,a(n)^2$ if $n$ is odd and $\tau(u)=6n\,a(n)^2$ if $n$ is even.

\subsection{Dihedral graph $\mathcal{D}_{n}=Cay(\mathbb{D}_{n}, b^{\pm1}, b^{\pm2},  a b,a b^3,a b^5).$}\label{example2}

\noindent$\textbf{1}^{\circ}.$ \textbf{The number of spanning trees.} By Theorem~\ref{theorem1} 
$$\tau(n)=\frac{n \,t\, \eta^n}{q}|(-2)^n-1|\cdot|(-1/2)^n-1|\cdot|(4-\sqrt{15})^n-1|\cdot|(4+\sqrt{15})^n-1|,$$ 
where $t=3,\,\eta=2,\,q=54.$ Equivalently, by Theorem~\ref{theorem2} we get
$$\tau(n)=\frac{n\,2^{n}}{18}|2T(n,-\frac{5}{4})-2|\cdot|2T(n,4)-2|.$$

\noindent$\textbf{2}^{\circ}.$ \textbf{The asymptotics of $\tau(n).$} By Theorem~\ref{theorem3}, we have $\tau(n)\cong \frac{n}{18} A^n,$ where $A=4(4+\sqrt{15}).$\bigskip

\noindent$\textbf{3}^{\circ}.$ \textbf{The generating function of $\tau(n).$} From Theorem~\ref{theoremR1}, we get
$$F(x)=\sum\limits_{n=1}^\infty\tau(n)x^n$$

$$=\frac{6 (-1745300 + 4540750 u - 3003815 u^2 + 346990 u^3 + 171265 u^4 - 47660 u^5 + 4840 u^6 - 272 u^7 + 16 u^8)}{(2 + u)(8 + u)^2 (-5 + 2 u)^2 (265 - 80 u + 4 u^2)^2},$$  
where $u=\frac{1}{2}(2 x+\frac{1} {2 x} ).$
     

     \bigskip     

\noindent$\textbf{4}^{\circ}.$ \textbf{Divisibility by squares.} To see  the divisibility by squares we consider a few terms of   generating function
$$F(x)=3 x + 60 x^2 + 6561 x^3 + 192000 x^4 + 9149415 x^5 + 315059220 x^6+\cdots .$$
By Theorem \ref{lorenzini}, we have $\xi=(2\beta_{odd}+\gamma_{odd})(2\beta_{odd}+\gamma_{even})=(2\cdot1+3)(2\cdot1+0)=10.$ Hence, $\delta=10.$ So that there exists an integer sequence  $a(n)$ such that $\tau(u)=3n\,a(n)^2$ if $n$ is odd and $\tau(u)=30n\,a(n)^2$ if $n$ is even.

\section*{ACKNOWLEDGMENTS}
The first author was supported by NSFC, no. 11831004 and Shanghai Science and Technology Program [Project No. 22JC1400100]. The second and the third authors were supported by Mathematical Center in Akademgorodok under agreement No. 075-15-2019-1613 with the Ministry of Science and Higher Education of the Russian Federation. The fourth author is supported by  NSFC, no.12101125, Natural Science Foundation of Fujian Province [Project No.2021J05035].

\end{document}